\documentclass[a4paper,12pt]{article}
\usepackage{amscd}
\usepackage{amsmath,amsfonts,amssymb,amscd}
\usepackage{indentfirst,graphicx,epsfig}
\usepackage{graphicx,psfrag}
\usepackage{ifpdf}
\usepackage{float}

\input{epsf}

\newtheorem{thm}{Theorem}

\newtheorem{cor}{Corollary}

\def\qed{\hfill \nopagebreak\rule{5pt}{8pt}}

\begin{document}
\rule{0cm}{1cm}
\begin{center}
{\Large\bf Rainbow connection number, bridges and
radius\footnote{Supported by NSFC No.11071130.}}
\end{center}

\begin{center}
Jiuying Dong, Xueliang Li\\
Center for Combinatorics, LPMC\\
Nankai University, Tianjin 300071, P. R. China\\
Email:  jiuyingdong@126.com, lxl@nankai.edu.cn,
\end{center}

\begin{abstract}

Let $G$ be a connected graph. The notion \emph{the rainbow
connection number $rc(G)$} of a graph $G$ was introduced recently by
Chartrand et al. Basavaraju et al. showed that for every bridgeless
graph $G$ with radius $r$, $rc(G)\leq r(r+2)$, and the bound is
tight. In this paper, we prove that if $G$ is a connected graph, and
$D^{k}$ is a connected $k$-step dominating set of $G$, then $G$ has
a connected $(k-1)$-step dominating set $D^{k-1}\supset D^{k}$ such
that $rc(G[D^{k-1}])\leq rc(G[D^{k}])+\max\{2k+1,b_k\}$, where $b_k$
is the number of bridges in $ E(D^{k}, N(D^{k}))$. Furthermore, for
a connected graph $G$ with radius $r$, let $u$ be the center of $G$,
and $D^{r}=\{u\}$. Then $G$ has $r-1$ connected dominating sets $
D^{r-1}, D^{r-2},\cdots, D^{1}$ satisfying $D^{r}\subset
D^{r-1}\subset D^{r-2} \cdots\subset D^{1}\subset D^{0}=V(G)$, and
$rc(G)\leq \sum_{i=1}^{r}\max\{2i+1,b_i\}$, where $b_i$ is the
number of bridges in $ E(D^{i}, N(D^{i})), 1\leq i \leq r$. From the
result, we can get that if for all $1\leq i\leq r, b_i\leq 2i+1$,
then $rc(G)\leq \sum_{i=1}^{r}(2i+1)= r(r+2)$; if for all $1\leq
i\leq r, b_i> 2i+1$, then $rc(G)= \sum_{i=1}^{r}b_i$, the number of
bridges of $G$. This generalizes the result of
Basavaraju et al.\\[3mm]

\end{abstract}

\section{Introduction}

All graphs considered in this paper are simple, finite and
undirected. Undefined terminology and notations can be found in
\cite{Bondy-Murty}. Let $G$ be a graph, and $c: E(G) \rightarrow
\{1,2,\cdots,k\}, k\in N$ be an edge-coloring, where adjacent edges
may be colored the same. A graph $G$ is \emph{rainbow connected} if
for any pair of distinct vertices $u$ and $v$ of $G$, $G$ has a
$u-v$ path whose edges are colored with distinct colors. The minimum
number of colors required to make $G$ rainbow connected is called
its \emph{rainbow connection number}, denoted by $rc(G)$. These
concepts were introduced by Chartrand et al. in
\cite{Chartrand-Johns}, where they determined the rainbow connection
numbers of wheels, complete graphs and all complete multipartite
graphs. Many results involving some graph parameters were obtained.
Results involving the minimum degree were obtained in
\cite{Caro-Lev, Schiermeyer1, Krivelevich-Yuster, Chandran-Das}.
Results involving the parameters $\sigma_2$ and $\sigma_k(G)$ were
obtained in \cite{ Schiermeyer, Dong-Li}. In
\cite{Basavaraju-Chandran}, Basavaraju et al. showed that for every
bridgeless graph $G$ with radius $r$, $rc(G)\leq r(r+2)$, and the
bound is tight. As one can see, they did not consider graphs with
bridges. In this paper, we will consider graphs with bridges, and
$rc(G)$ is bounded by the number of bridges and radius of the
graphs. The following are our main results.

\begin{thm}
If $G$ is a  connected graph, and $D^{k}$ is a connected  $k$-step
dominating set of $G$, then $G$ has a connected $(k-1)$-step
dominating set $D^{k-1}\supset D^{k}$ such that $rc(G[D^{k-1}])\leq
rc(G[D^{k}])+\max\{2k+1,b_k\}$, where $b_k$ is the number of bridges
of $G$ in $ E(D^{k}, N(D^{k}))$.
\end{thm}

\begin{thm}
For a connected graph $G$ with radius $r$, let $u$ be the center of
$G$, and $D^{r}=\{u\}$. Then $G$ has $r-1$ connected dominating sets
$ D^{r-1}, D^{r-2},\cdots, D^{1}$ satisfying $D^{r}\subset
D^{r-1}\subset D^{r-2}\cdots \subset D^{1}\subset D^{0}=V(G)$, and
$rc(G)\leq \sum_{i=1}^{r}\max\{2i+1,b_i\}$, where $b_i$ is the
number of bridges in $ E(D^{i}, N(D^{i})), 1\leq i \leq r$.
\end{thm}

Note that if for all $1\leq i\leq r, b_i\leq 2i+1$, then $rc(G)\leq
\sum_{i=1}^{r}(2i+1)= r(r+2)$; if for all $1\leq i\leq r, b_i>
2i+1$, then $rc(G)= \sum_{i=1}^{r}b_i$, the number of bridges of
$G$. This generalizes the result of Basavaraju et al.

\section{Preliminaries}

For two subsets $X$ and $Y$ of $V$, an $(X,Y)$-path is a path which
connects a vertex of $X$ and a vertex of $Y$, and whose internal
vertices belong to neither $X$ nor $Y$. We use $E[X,Y]$ to denote
the set of edges of $G$ with one end in $X$ and the other end in
$Y$, and $e(X,Y)=|E[X,Y]|$.

Let $G$ be a connected graph. The eccentricity of a vertex $v$ is
$ecc(v)=\max_{x\in V(G)}d_G(v,x)$. The radius of $G$ is
$rad(G)=\min_{x\in V(G)}ecc(x)$. The diameter of $G$ is $diam(G)=
\max_{x\in V(G)}ecc(x)$. Let $S\subseteq V(G)$. The $k$-step open
neighborhood of $S$ is $N^{k}(S)=\{v\in V(G)|d(v,S)=k, k\in Z, k\geq
0\}$. Generally speaking, $N^{1}(S)=N(S), N^{0}(S)=S,
N^{k}[S]=N^{k}(S)\cup S$. If every vertex in $G$ is at a distance at
most $k$ from $S$, we say that $S$ is a $k$-step dominating set. If
$S$ is connected, then $S$ is a connected $k$-step dominating set.

The following definitions are needed in our proof. Let $D^{k}$ be a
connected $k$-step dominating set. A $D^{k}$-ear is a path
$P=v_0v_1\cdots v_p$ in $G$ such that $P\cap D^{k}=\{v_0,v_p\}$.
When $v_0=v_p$, $P$ is a closed $D^{k}$-ear. Moreover, we say that
$P$ is an eager $D^{k}$-ear, if $P$ is a shortest $D^{k}$-ear
containing $v_0v_1$. Given $2k+1$ distinct colors, for convenience,
we denote them by $1,2,3,\cdots,2k+1$. We say that $P$ is evenly
colored, if either the edges of $P$ are colored in this way:
$c(v_0v_1)=1, c(v_1v_2)=2, c(v_2v_3)=3,\cdots,
c(v_{\lceil\frac{p}{2}\rceil-1}v_{\lceil\frac{p}{2}\rceil})=
\lceil\frac{p}{2}\rceil, c(v_{\lceil\frac{p}{2}\rceil}
v_{\lceil\frac{p}{2} \rceil+1})=2k+2-\lfloor\frac{p}{2}\rfloor,
c(v_{\lceil\frac{p}{2}\rceil+1}v_{\lceil\frac{p}{2}\rceil+2})
=2k+3-\lfloor\frac{p}{2}\rfloor, \cdots, c(v_{p-2}v_{p-1})=2k,
c(v_{p-1}v_p)=2k+1$, or the edges of $P$ are colored in another way:
$c(v_0v_1)=2k+1, c(v_1v_2)=2k, c(v_2v_3)=2k-1,\cdots,
c(v_{p-\lceil\frac{p}{2}\rceil-2}v_{p-\lceil\frac{p}{2}\rceil-1})
=2k+2-\lfloor\frac{p}{2}\rfloor,
c(v_{p-\lceil\frac{p}{2}\rceil-1}v_{p-\lceil\frac{p}{2}\rceil})
=\lceil\frac{p}{2}\rceil, \cdots, c(v_{p-2}v_{p-1})=2,
c(v_{p-1}v_p)=1$. In the proofs later, for convenience, we say that
$P$ is evenly colored, if either the edges of $P$ are colored
$1,2,\cdots, \lceil\frac{p}{2}\rceil,2k+2-\lfloor\frac{p}{2}\rfloor,
2k+1-\lfloor\frac{p}{2}\rfloor,\cdots, 2k,2k+1$ in this order, or
the edges of $P$ are colored $2k+1,2k,\cdots,
2k+2-\lfloor\frac{p}{2}\rfloor, \lceil\frac{p}{2}\rceil,
\lceil\frac{p}{2}\rceil-1,\cdots, 3,2,1$ in this order.

\section{The proofs of our theorems}

\parindent 0pt

\noindent{\bf The proof of Theorem 1}:

If $G$ is a tree, then each edge of $G$ is bridge, The result is
obvious. Hence we may assume that $G$ is not a tree.

The following, we let $D^{k}$ be a connected $k$-step dominating set
of $G$. Then $G$ has $k$ mutually disjoint subsets $N^{1}(D^{k}),
N^{2}(D^{k}),\cdots,N^{k}(D^{k})$ and
$V(G)=\bigcup_{i=0}^{k}N^{i}(D^{k})$.

\textbf{Claim 1:} If $\exists x\in N(D^{k}), y\in D^{k}$ such that
$xy$ is bridge, then we have $d_{G[N[D^{k}]]}(x)=1$.

If $\exists y'\in D^{k}, y'\neq y$, such that $xy'\in E(G)$. As
$G[D^{k}]$ is connected, $G[D^{k}]$ has a path connecting $y,y'$.
Hence $xy$ is in a cycle, a contradiction to $xy$ being a bridge. If
$\exists x_1\in N( D^{k})$ such that $xx_1\in E(G)$, as there exists
some vertex $y_1\in D^{k}$ satisfying $x_1y_1\in E(G)$ ($y_1$ may be
$y$), then $yxx_1y_1$ is a path, and $G[D^{k}]$ has a path
connecting $y,y_1$, that is, $xy$ is in some cycle, a contradiction.
Hence $d_{G[N[D^{k}]]}(x)=1$.

Let $x_1y_1, x_2y_2,\cdots,x_{b_r}y_{b_r}$ be all the bridges in $E(
N( D^{k}),  D^{k})$, where $x_i\in  N( D^{k}), y_i\in D^{k}, 1\leq
i\leq b_r $. Set $B=\{x_1,x_2,\cdots, x_{b_r}\}$,
$B_{E}=\{x_1y_1,x_2y_2, \cdots, x_{b_r}y_{b_r}\}$, $D_1= D^{k}\cup
B$.

Let  $D^{k}$ be a connected $k$-step dominating set,
we rainbow color $G[D^{k}]$ with $rc(G[D^{k}])$ colors.
 If $N^( D^{k})=B$, then $D_1$ is a connected $(k-1)$-
step dominating set. Set $ D^{k-1}=D_1$, we use $b_r$ fresh colors
to color these $b_r$ bridges, respectively. Hence
$rc(G[D^{k-1}])\leq rc(G[D^{k}])+b_r$, and the theorem follows.

So we may assume $N( D^{k})\setminus B\neq \emptyset$. For any
vertex $v_1\in N( D^{k})\setminus B$, and any edge $v_0v_1\in E(v_1,
D^{k})$, $v_0\in D^{k}$, as $v_0v_1$ is not a bridge, $v_0v_1$ is in
some cycle. Hence we may let $P=v_0v_1v_2\cdots
v_tv_{t+1}v_{t+2}\cdots v_{t+m}v_{t+m+1} \cdots v_{p-1}v_p$ be an
eager $D^{k}$-ear.

\textbf{Claim 2:} $|P|\leq 2k+1$.

It mainly depends on the following Claim 2.1 and Claim 2.2.

\textbf{Claim 2.1:}  If $v_t,v_{t+1}\in N^{t}( D^{k}), v_i\in N^{i}(
D^{k}), 0\leq i\leq t-1, t\geq 1$, then $v_{t+2}\in N^{t-1} (
D^{k})$, and $P$ does not have two vertices $v_{t+m},v_{t+m+1}$ in
some $N^{j}( D^{k})$, where $m\geq 2, 1\leq j\leq t-1$.

If $v_1,v_2\in N( D^{k}), v_1v_2\in E(G)$, then there
exists $v_3\in D^{k}$ ($v_3$ can be $v_0$) such that
$v_0v_1v_2v_3$ is an eager $D^{k}$-ear.

So we may assume $t\geq 2$. Suppose that, to the contrary,
$v_{t+2}\in N^{t}( D^{k})$ or $v_{t+2}\in  N^{t+1}( D^{k})$. If
$v_{t+1}v_{t-1}\in E(G)$, then we replace $P$ by a shorter path
$P'=v_0v_1v_2\cdots v_{t-2}v_{t-1}v_{t+1} v_{t+2}\cdots $\\
$v_{t+m}v_{t+m+1}\cdots v_{p-1}v_p$. If $v_{t+1}v_b\in E(G), b\in
N^{t-1}( D^{k})\cap (P\setminus \{v_{t-1}\})$, then we replace $P$
by a shorter path $P'=v_0v_1v_2\cdots
v_{t-1}v_tv_{t+1}v_bv_{b+1}\cdots v_{t+m}v_{t+m+1}\cdots
v_{p-1}v_p$, a contradiction to $P$ being an eager $D^{k}$-ear.
Hence $v_{t+2}\in N^{t-1}( D^{k})$.

Suppose that $P$ has two vertices $v_{t+m},v_{t+m+1}$ in some
$N^{j}( D^{k})$ where $m\geq 2, 1\leq j\leq t-1$. If $j=1$, that is
$v_{t+m},v_{t+m+1}\in N( D^{k})$, $v_{t+m}v_{t+m+1}\in E(G)$,
because there is some vertex $v_{p_1}\in  D^{k}$($v_{p_1}$ may be
$v_0$) such that $v_{t+m}v_{p_1}\in E(G)$, then we replace $P$ by a
shorter path $P'=v_0v_1v_2\cdots v_{t-1}v_tv_{t+1}\cdots
v_{t+m-1}v_{t+m} v_{p_1}$, a contradiction. So we may assume $2\leq
j\leq t-1$. If $v_{t+m}v_{j-1}\in E(G)$, then we replace $P$ by a
shorter path $P'=v_0v_1v_2 \cdots v_{j-1}v_{t+m}v_{t+m+1}\cdots
v_{p-1}v_p$. If $v_{t+m}v_a\in E(G)$, where $a\in N^{j-1}( D^{k})
\cap (P\setminus \{v_{j-1}\})$, then we replace $P$ by a shorter
path $P'=v_0v_1v_2\cdots v_{t-1}v_tv_{t+1}\cdots
v_{t+m}v_av_{a+1}\cdots v_{p-1}v_p$, a contradiction to $P$ being an
eager $D^{k}$-ear. Hence $P$ does not have two vertices
$v_{t+m},v_{t+m+1}$ in some $N^{j}( D^{k})$ where $m\geq 2, 1\leq
j\leq t-1$. Claim 2.1 is true.\qed

\textbf{Claim 2.2:} If $v_i\in N^{i}( D^{k}), 0\leq i\leq t, t\geq
2$ and $v_{t+1}\in N^{t-1}( D^{k})$, then $P$ does not have two
vertices $v_{t+m},v_{t+m+1}$ in some $N^{j}( D^{k})$ where $m\geq 1,
1\leq j\leq t-1$.

Suppose that, to the contrary, $P$ has two vertices
$v_{t+m},v_{t+m+1}$ in some $N^{j}( D^{k})$ where $m\geq 1, 1\leq
j\leq t-1$. The proof of Claim 2.2 is similar to the proof in the
latter part of Claim 2.1, and so Claim 2.2 is also true. \qed

By Claim 2.1 and Claim 2.2, we can get $|P|\leq 2k+1$, in which
equality holds if and only if $t=k$ and $v_kv_{k+1}\in N^{k}(D^{k}),
v_{k+2}\in N^{k-1}(D^{k})$, and so Claim 2 is true. \qed

In the following we will construct a connected $(k-1)$-step
dominating set $D^{k-1}$ such that $G[D^{k-1}]$ is rainbow
connected. Since $E(D^{k}, N( D^{k})\setminus B)$ has no bridges,
for each edge $e$ of $E(D^{k}, N( D^{k})\setminus B)$, $e$ must be
in some cycle, and so there exists an eager $D^{k}$-ear $P$
containing $e$. Thus we may construct a sequence of sets $D_1\subset
D_2\subset D_3\subset\cdots\subset D_t=D^{k-1}$, where $D_2=D_1\cup
P_1$, $D_3=D_2\cup P_2,\cdots, D_t=D_{t-1}\cup P_{t-1}$,
$P_1,P_2,\cdots,P_t$ are all eager $D^{k}$-ears. We color the new
edges in every induced graph $G[D_i]$ such that every $x\in
D_i\setminus D_1$ lies in an evenly colored eager $D^{k}$-ear in
$G[D_i]$ for all $1\leq i \leq t$.

For $i=1$, it is obvious. If for some $D_i$, $N(D^{k})\subset D_i$,
note that for $1\leq j\leq i-1$, $ N( D^{k})\not\subset D_j$, then
$D_i$ is a connected $(k-1)$-step dominating set. We stop the
procedure and set $D^{k-1}=D_i$, and evenly color the edges of
$P_{i-1}$, and color the remaining uncolored new edges of $G[D_i]$
with the used colors. Otherwise, we will construct $D_{i+1}$ as
follows:

We choose any edge $x_0x_1\in E(D^{k}, N( D^{k})\setminus D_i)$,
$x_0\in D^{k}, x_1\in N^{1}( D^{k})\setminus D_i$. If $P$ is an
eager $D^{k}$-ear containing $x_0x_1$, and $P\cap(D_i\setminus
D_1)=\emptyset$, then we set $D_{i+1}=D_i\cup P$, and evenly color
$P$, for the uncolored new edges of $G[D_{i+1}]$, we color them
randomly with the used colors. Otherwise, the eager $D^{k}$-ear $P$
containing $x_0x_1$ must satisfy $P\cap(D_i\setminus
D_1)\neq\emptyset$. Assume $P_1\subset P$, and let $P_1=x_0x_1\cdots
x_l$, $P_1\cap(D_i\setminus D_1)=\{x_l\}$. As $x_l\in D_i\setminus
D_1$, $x_l$ is in an evenly colored eager $D^{k}$-ear $Q$. Let $Q_1$
be the shorter segment of $Q$ respect to $x_l$. Then $P=P_1\cup Q_1$
is the eager $D^{k}$-ear containing $x_0x_1$. We know that $Q$ is
evenly colored. If $Q_1$ is colored by the colors from
$\{2k+1,2k,2k-1,\cdots, 2k+2-\lfloor\frac{|Q|}{2}\rfloor\}$, then we
will evenly color $P$ by $1,2,3,\cdots \lceil\frac{|P|}{2}\rceil,
2k+2-\lfloor\frac{|P|}{2}\rfloor,\cdots, 2k,2k+1$ in that order,
here $c(x_0x_1)=1$. If $Q_1$ is colored by the colors from
$\{1,2,3,\cdots, \lceil\frac{|Q|}{2}\rceil\}$, then we will evenly
color $P$ by $2k+1,2k, 2k-1,\cdots,
2k+2-\lfloor\frac{|P|}{2}\rfloor,
\lceil\frac{|P|}{2}\rceil,\cdots,3,2,1$ in that order, here
$c(x_0x_1)=2k+1$. Hence $P$ is evenly colored. Set $D_{i+1}=D_i\cup
P$. For the uncolored new edges of $G[D_{i+1}]$, we color them
randomly with the used colors. Clearly, every $x\in D_{i+1}\setminus
D_1$ lies in an evenly colored eager $D^{k}$-ear in $G[D_{i+1}]$.

Thus, we have constructed a connected $(k-1)$-step dominating set
$D^{k-1}$, and every edge of $G[D^{k-1}\setminus B]$ is colored.

\textbf{Claim 3:} $G[D^{k-1}\setminus D^{k}]$ has no bridges.

Suppose that $xy\in G[D^{k-1}\setminus D^{k}]$ is a bridge.

By Claim 1, we know that if $x\in B$, then $y\not\in B$, and if
$y\in B$, then $x\not\in B$. Hence we will consider the following
two cases: If $x$ is in some eager $D^{k}$-ear $P$, $y$ is in some
eager $D^{k}$-ear $Q$ ($P$ can be $Q$), then besides $xy$, there is
still another path connecting $x$ and $y$, so $xy$ is in a cycle. If
$x\in B$, $y$ is in some eager $D^{k}$-ear $Q$, then $xy$ is also in
some cycle, a contradiction. \qed

Now, we are ready for coloring $B_E$: If $b_k\leq 2k+1$, then we use
$b_k$ different colors from $\{1,2,\cdots,2k+1\}$ to color each edge
of $B_E$, respectively. If $b_k >2k+1$, then we first use colors
$1,2,\cdots,2k+1$ to color any $2k+1$ edges of $B_E$, respectively,
then we use $b_k-(2k+1)$ fresh colors to color the remaining
uncolored edges, respectively.

In the following we claim that $G[D^{k-1}]$ is rainbow connected.
For any two vertices $x,y\in D_1$, we know that $x,y$ is rainbow
connected. For $x\in D^{k-1}\setminus D_1, y\in D^{k}$, as $x$ is in
an eager $D^{k}$-ear $P$, let $P\cap D^{k}=y_1$. In $D^{k}$, there
exists a rainbow path connecting $y,y_1$. For $x\in D^{k-1}\setminus
D_1, y\in B$, we know that $x$ is in an evenly eager $D^{k}$ ear
$P$. If the bridge $yy_1\in B_E(y_1\in D^{k}$) is colored by $c_y$
which is also in $P$, then we choose the segment (which does not
contain the color $c_y$) connecting $x$ to $D^{k}$. If the bridge
$yy_1\in B_E,(y_1\in D^{k}$) is colored by $c_y$ which is not in
$P$, then we arbitrarily choose a segment of $P$ connecting $x$ to
$D^{k}$, we can also find a $x-y$ rainbow path.

For $x\in D^{k-1}\setminus D_1, y\in D^{k-1}\setminus D_1$, since
$x$ and $y$ are both in evenly colored eager $D^{k}$-ears, let $x\in
P,y\in Q$, $P, Q$ are evenly colored eager $D^{k}$-ears. If $P=Q$,
then $x, y$ is rainbow connected. Hence we may assume $P\neq Q$. Let
$P=x_0x_1\cdots x_i(x)x_{i+1}\cdots x_p$, $Q=y_0y_1\cdots
y_j(y)y_{j+1}\cdots y_q$. We distinguish two cases to show that
$x,y$ is rainbow connected.

\textbf{Case 1:} $P$ and $Q$ are internally disjoint.

We assume that $x_0x_1,\cdots,x_{\lceil\frac{p}{2}\rceil}$ and
$y_0y_1,\cdots,y_{\lceil\frac{q}{2}\rceil}$ are colored by the
colors from $\{1,2,3,\cdots,\\k+1\}$, respectively. The other three
coloring cases can be discussed in a similar way. We distinguish
four subcases to demonstrate that there is an $x-y$ rainbow path.

\textbf{Subcase 1.1:}  $i\leq\lfloor\frac{p}{2}\rfloor$,
$j>\lfloor\frac{q}{2}\rfloor$.

We join $x=x_ix_{i-1}\cdots x_0$ to the $x_0-y_q$ rainbow path in
$G[D^{k}]$ followed by $y_q y_{q-1}\cdots y_j=y$. As the edges of
$x=x_ix_{i-1}\cdots x_0$ are colored by the colors from
$\{1,2,\cdots, k+1\}$, the edges of $y_q y_{q-1}\cdots y_j=y$ are
colored by the colors from $\{2k+1,2k, \cdots, k+2\}$. Hence it is
an $x-y$ rainbow path.

\textbf{Subcase 1.2:}  $i> \lfloor\frac{p}{2}\rfloor$,
$j\leq \lfloor\frac{q}{2}\rfloor$.

We join $y=y_jy_{j-1}\cdots y_0$ to the $y_0-x_p$ rainbow path in
$G[D^{k}]$ followed by $x_p x_{p-1}\cdots x_i=x$. It is also an
$x-y$ rainbow path.

\textbf{Subcase 1.3:}  $i\leq \lfloor\frac{p}{2}\rfloor$,
$j\leq \lfloor\frac{q}{2}\rfloor$.

If $i< j$, we join $x=x_ix_{i-1}\cdots x_0$ to the $x_0-y_q$ rainbow
path in $G[D^{k}]$ followed by $y_qy_{q-1}\cdots y_j=y$. As the
edges of $x=x_ix_{i-1}\cdots x_0$ are colored by $i,i-1,\cdots, 1$,
the edges of $y_q y_{q-1}\cdots y_j=y$ are colored by the colors
$2k+1,2k, \cdots, j$. It is an $x-y$ rainbow path. If $i\geq j$, we
join $y=y_jy_{j-1}\cdots y_0$ to the $y_0-x_p$ rainbow path in
$G[D^{k}]$ followed by $x_p x_{p-1}\cdots x_i=x$. As the edges of
$y=y_jy_{j-1}\cdots y_0$ are colored by the colors $\{j,j-1,\cdots,
1\}$, the edges of $x_px_{p-1}\cdots x_i=x$ are colored by the
colors $\{2k+1,2k, \cdots, i\}$, it is also an $x-y$ rainbow path.

\textbf{Subcase 1.4:}  $i > \lfloor\frac{p}{2}\rfloor$,
$j > \lfloor\frac{q}{2}\rfloor$.

If $p-i\leq q-j$, then we join $x=x_ix_{i+1}\cdots x_p$ to the
$x_p-y_0$ rainbow path in $G[D^{k}]$ followed by $y_0y_1,
\cdots,y_j=y$. If $p-i> q-j$, we join $y=y_jy_{j+1}\cdots y_q$ to
the $y_q-x_0$ rainbow path in $G[D^{k}]$ followed by $x_0x_1\cdots
x_i=x$. So we find an $x-y$ rainbow path.

\textbf{Case 2:} $P$ and $Q$ are internally joint.

According to the construction and the coloring of $D^{k-1}$, we may
assume that $P\subset D_{i_1}$, $Q\subset D_{i_2}$, and $i_1>i_2$,
$x_l$ is the first internal vertex of $P$ in $Q$. If
$x_px_{p-1}\cdots x_{l+1}x_l=y_qy_{q-1}\cdots y_{l+1}y_l$, then the
case is similar to Case 1 in essence. So we may assume
$x_px_{p-1}\cdots x_{l+1}x_l= y_0y_1,\cdots,y_{p-l}$. We also
distinguish four subcases to show that there is an $x-y$ rainbow
path.

Without loss of generality, assume that the edges of $y_0y_1\cdots
y_{\lceil\frac{q}{2}\rceil}$ are colored by $1,2,
\cdots,{\lceil\frac{q}{2}\rceil}$. According to the coloring of
$D^{k-1}$, the edges of $x_px_{p-1}\cdots
x_{\lfloor\frac{p}{2}\rfloor}$ are also colored by the colors from
$\{1,2, \cdots, k+1\}$, and the edges of $x_0x_1\cdots
x_{\lceil\frac{p}{2}\rceil}$ are colored by the colors from
$\{2k+1,2k, \cdots, k+2\}$.

\textbf{Subcase 2.1:} $i\leq\lfloor\frac{p}{2}\rfloor$,
$j>\lfloor\frac{q}{2}\rfloor$.

If $i< q-j$, then we join $x=x_ix_{i-1}\cdots x_0$ to the $x_0-y_0$
rainbow path in $G[D^{k}]$ followed by $y_0y_1\cdots y_j=y$. If
$i\geq q-j$, then we join $y=y_jy_{j+1}\cdots y_q$ to the $y_q-x_p$
rainbow path in $G[D^{k}]$ followed by $x_px_{p-1}\cdots x_i=x$. We
find the required $x-y$ rainbow path.

\textbf{Subcase 2.2:}  $i> \lfloor\frac{p}{2}\rfloor$,
$j\leq \lfloor\frac{q}{2}\rfloor$.

If $p-i\leq j$, then we join $x=x_ix_{i+1}\cdots x_p$ to the
$x_p-y_q$ rainbow path in $G[D^{k}]$ with $y_qy_{q-1}\cdots y_j=y$.
If $p-i>j$, we join $y=y_jy_{j-1}\cdots y_0$ to the $y_0-x_0$
rainbow path in $G[D^{k}]$ followed by $x_0x_1\cdots x_i=x$. We also
find the required $x-y$ rainbow path.

\textbf{Subcase 2.3:}  $i\leq \lfloor\frac{p}{2}\rfloor$, $j\leq
\lfloor\frac{q}{2}\rfloor$.

we join $x=x_ix_{i-1}\cdots x_0$ to the $x_0-y_0$ rainbow path in
$G[D^{k}]$ followed by $y_0y_1\cdots y_j=y$.

\textbf{Subcase 2.4:}  $i > \lfloor\frac{p}{2}\rfloor$,
$j > \lfloor\frac{q}{2}\rfloor$.

We join $x=x_ix_{i+1}\cdots x_p$ to the $x_p-y_q$ rainbow path in
$G[D^{k}]$ followed by $y_qy_{q-1}\cdots y_j=y$.

Hence, for any two vertices $x,y\in D^{k-1}\setminus D_1$, there is
rainbow path connecting $x$ and $y$. Thus, we have constructed a
connected $D^{k-1}$ from $D^{k}$, and $rc(G[D^{k-1}])\leq
rc(G[D^{k}])+\max\{2k+1,b_k\}$.

Hitherto, the proof of Theorem 1 has been completed. \qed

\noindent{\bf The proof of Theorem 2}:

Let $u$ be the center of $G$, and set $D^{r}=\{u\}$. Then $D^{r}$ is
an $r$-step dominating set of $G$, and $rc(G[D^{r}])=0$. By making
use of Theorem 1, we may construct $D^{r-1}, D^{r-2},\cdots,
D^{2},D^{1}$ such that $D^{r}\subset D^{r-1}\subset D^{r-2}\cdots
\subset D^{1}\subset D^{0}=V(G)$, and we have
$$rc(G[D^{r-1}])\leq rc(G[D^{r}])+\max\{2r+1, b_r\}$$
$$rc(G[D^{r-2}])\leq rc(G[D^{r-1}])+\max\{2(r-1)+1, b_{r-1}\}$$
$$\cdots$$
$$rc(G[D^{0}])\leq rc(G[D^{1}])+\max\{2+1, b_1\}$$
where $rc(G[D^{0}])=rc(G)$, for $1\leq i \leq r$, $b_i$ is the
number of bridges in $E(D^{i}, N(D^{i}))$. Thus we get that
$rc(G)\leq rc(G[D^{r}])+\sum_{i=1}^{r}\max\{2i+1,b_i\}
=\sum_{i=1}^{r}\max\{2i+1,b_i\}$.

This completes the proof of Theorem 2.\qed

By Claim 3, the subgraph $G[D^{i-1}\setminus D^{i}]]$ has no
bridges. Hence we immediately obtain the following corollary.

\begin{cor}
The number of bridges of $G$ is equal to $\sum_{i=1}^{r} b_i$.
\end{cor}

\end{document}